\documentclass{article}[10pt]
\usepackage[T1]{fontenc}
\usepackage[latin1]{inputenc}
\usepackage[english]{babel}
\usepackage{hyperref}
\usepackage[dvips]{graphics}
\usepackage{amsmath}
\usepackage{amssymb}
\usepackage{amsopn}
\usepackage{amsbsy}
\usepackage{amstext}
\usepackage{latexsym}
\usepackage{amsfonts}
\usepackage{array}
\usepackage{epsfig}
\usepackage{mathrsfs}
\title{Upper bounds for the order of an additive basis obtained by removing a finite subset of a given basis}
\author{Bakir FARHI}
\date{}

\newtheorem{thm}{Theorem}[section]
\newtheorem{prop}[thm]{Proposition}
\newtheorem{lemma}[thm]{Lemma}
\newtheorem{rmq}[thm]{Remark}
\newtheorem{coll}[thm]{Corollary}

\let\epsilon=\varepsilon
\def\EMts{\mspace{.3mu}}
\def\diam{{\rm diam}}

\def\gcd{{\rm gcd}}
\def\G{{\rm G}}
\def\stab{{\rm stab}}
\def\d{{\underline{\bf{d}}}}
\def\nb#1{{\left\vert{\EMts\EMts #1 \EMts\EMts}\right\vert}}

\def\EMdash{\leavevmode\hbox to 7.5mm{\vrule height .63ex depth -.59ex
    width 5.4mm\hfill}}

\begin{document}
\maketitle \vspace{-6cm}
\begin{flushleft}
{\it J. Number Theory}, \\
{\bf 128} (2008), p. 2214-2230.
\end{flushleft}~\vspace{3cm}
\begin{center}
{\tt bakir.farhi@gmail.com}
\end{center}
\begin{abstract}
Let $A$ be an additive basis of order $h$ and $X$ be a finite
nonempty subset of $A$ such that the set $A \setminus X$ is still a
basis. In this article, we give several upper bounds for the order
of $A \setminus X$ in function of the order $h$ of $A$ and some
parameters related to $X$ and $A$. If the parameter in question is
the cardinality of $X$, Nathanson and Nash already obtained some
of such upper bounds, which can be seen as polynomials in $h$ with
degree $(\nb{X} + 1)$. Here, by taking instead of the cardinality
of $X$ the parameter defined by $d := \frac{\diam(X)}{\gcd\{x - y
~|~ x , y \in X\}}$, we show that the order of $A \setminus X$ is
bounded above by $(\frac{h (h + 3)}{2} + d \frac{h (h - 1) (h +
4)}{6})$. As a consequence, we deduce that if $X$ is an arithmetic
progression of length $\geq 3$, then the upper bounds of Nathanson
and Nash are considerably improved. Further, by considering more
complex parameters related to both $X$ and $A$, we get upper
bounds which are polynomials in $h$ with degree only $2$.
\end{abstract}
\vspace{1mm}
\section{Introduction}~

An additive basis (or simply a basis) is a subset $A$ of
$\mathbb{Z}$, having a finite intersection with $\mathbb{Z}^-$ and
for which there exists a natural number $h$ such that any
sufficiently large positive integer can be written as a sum of $h$
elements of $A$. The smaller number $h$ satisfying this property
is called ``the order of the basis $A$'' and we note it $\G(A)$.
If $A$ is a basis of order $h$ and $X$ is a finite nonempty subset
of $A$ such that $A \setminus X$ is still a basis, the problem
dealt with here is to find upper bounds for the order of $A
\setminus X$ in function of the order $h$ of $A$ and parameters
related to $X$ (resp. $X$ and $A$). The particular case when $X$
contains only one element, say $X = \{x\}$, was studied for the
first time by Erd\"os and Graham \cite{eg}. These two last authors
showed that $\G(A \setminus \{x\}) \leq \frac{5}{4} h^2 +
\frac{1}{2} h \log{h} + 2 h$. After hem, several works followed in
order to improve this estimate: In his Thesis, by using Kneser's
theorem (see e.g. \cite{k} or \cite{hr}), Grekos \cite{g} improved
the previous estimate to $\G(A \setminus \{x\}) \leq h^2 + h$. By
still using Kneser's theorem but in a more judicious way, Nash
\cite{n} improved the estimate of Grekos to $\G(A \setminus \{x\})
\leq \frac{1}{2} (h^2 + 3 h)$. Finally, by combining Kneser's
theorem with some new additive methods, Plagne \cite{p} obtained
the refined estimate $\G(A \setminus \{x\}) \leq \frac{h (h +
1)}{2} + \lceil\frac{h - 1}{3}\rceil$, which is best known till
now. Plagne conjectured that $\G(A \setminus \{x\}) \leq \frac{h
(h + 1)}{2} + 1$, but this has not yet been proved. Notice also
that the optimality of such estimates was discussed by different
authors (see e.g. \cite{eg}, \cite{g}, \cite{g2}, \cite{p}).

The general case of the problem was studied by Nathanson and Nash
(see e.g. \cite{nat}, \cite{n2}, \cite{nn} and \cite{n}). For $h ,
k \in \mathbb{N}$, these two authors noted $G_k(h)$ the maximum of
all the natural numbers $\G(A \setminus X)$, where $A$ is an
additive basis of order $h$ and $X$ is a subset of $A$ with
cardinality $k$ such that $A \setminus X$ is still a basis. In
\cite{nn}, they proved that $G_k(h)$ has order of magnitude $h^{k
+ 1}$. Indeed, they showed that
$$\left(\frac{h}{k + 1}\right)^{k + 1} + O(h^k) \leq G_k(h) \leq \frac{2}{k!}
h^{k + 1} + O(h^k)$$ (see Theorem 4 of \cite{nn}).\\
Since then, the above bounds of $G_k(h)$ were improved. In
\cite{x}, Xing-de Jia showed that
$$G_k(h) \geq \frac{4}{3}\left(\frac{h}{k + 1}\right)^{k + 1} + O(h^k)$$
and in \cite{n}, Nash obtained the following
\begin{thm}[\cite{n}, Proposition 3 simplified]\label{t1}
Let $A$ be a basis and $X$ be a finite subset of $A$ such that $A
\setminus X$ is still a basis. Then, noting $h$ the order of $A$
and $k$ the cardinality of $X$, we have:
$$\G(A \setminus X) \leq (h + 1) \binom{h + k - 1}{k} - k \binom{h + k - 1}{k + 1} .$$
\end{thm}
Actually, the original estimate of Nash (Proposition 3 of
\cite{n}) is that $\G(A \setminus X) \leq \binom{h + k - 1}{k} +
\sum_{i = 0}^{h - 1}\binom{k + i - 1}{i}(h - i)$. But we can
simplify this by remarking that for all $i \in \mathbb{N}$, we
have:
$$\binom{k + i - 1}{i} = \binom{k + i}{i} - \binom{k + i - 1}{i - 1}$$
and
$$i\binom{k + i - 1}{i} = k\binom{k + i - 1}{i - 1} = k\left\{\binom{k + i}{i - 1} - \binom{k + i - 1}{i - 2}\right\} .$$
Consequently, we have:
$$
\sum_{i = 0}^{h - 1}\binom{k + i - 1}{i}(h - i) ~=~ h\sum_{i =
0}^{h - 1} \binom{k + i - 1}{i} - \sum_{i = 0}^{h - 1} i\binom{k +
i - 1}{i}~~~~~~~~~~~~~~~~~~~~~~~~~~~~
$$
\begin{eqnarray*}
& = & h \sum_{i = 0}^{h - 1}\left\{\binom{k + i}{i} - \binom{k + i
- 1}{i - 1}\right\} - k \sum_{i = 0}^{h - 1}\left\{\binom{k + i}{i
- 1} - \binom{k + i - 1}{i - 2}\right\} \\
& = & h \binom{h + k - 1}{h - 1} - k \binom{h + k - 1}{h - 2} \\
& = & h \binom{h + k - 1}{k} - k \binom{h + k - 1}{k + 1} ,
\end{eqnarray*}
which leads to the estimate of Theorem \ref{t1}.

In Theorem \ref{t1}, the upper bound of $\G(A \setminus X)$ is
easily seen to be a polynomial in $h$ with leading term
$\frac{h^{k + 1}}{(k + 1)!}$, thus with degree $(k + 1)$. In this
paper, we show that it is even possible to bound from above $\G(A
\setminus X)$ by a polynomial in $h$ with degree constant ($3$ or
$2$) but with coefficients depend on a new parameter other the
cardinality of $X$. By setting
$$d := \frac{\diam(X)}{\delta(X)} ,$$ where $\diam(X)$ denotes the
usual diameter of $X$ and $\delta(X) := \gcd\{x - y ~|~ x , y \in
X\}$, we show that
$$\G(A \setminus X) \leq \frac{h (h + 3)}{2} + d \frac{h (h - 1) (h + 4)}{6} ~~~~~~ \text{(see Theorem \ref{t3}).}$$
Also, by setting
$$\eta := \min_{\begin{array}{c}\scriptstyle{a , b \in A \setminus X , a \neq b} \\ \scriptstyle{\nb{a - b} \geq
\diam(X)}\end{array}} \nb{a - b} ,$$ we show that
$$\G(A \setminus X) \leq \eta (h^2 - 1) + h + 1 ~~~~~~\text{(see Theorem \ref{t4}).}$$
Finally, by setting
$$\mu := \min_{y \in A \setminus X} \diam(X \cup \{y\}) ,$$
we show that $$\G(A \setminus X) \leq \frac{h \mu (h \mu + 3)}{2}
~~~~~~\text{(see Theorem \ref{t5}).}$$
 It must be noted that this last estimate is obtained by an
 elementary way as a consequence of Nash' theorem while the
 two first estimates are obtained by applying Kneser's theorem
 with some differences with \cite{n}.

 In practice, when $h$ and $k$ are large enough, it often happens
 that our estimates are better than that of Theorem \ref{t1}. The more interesting corollary is when $X$ is an
 arithmetic progression: in this case we have $d = k - 1$, implying from our first estimate an
 improvement of Theorem \ref{t1}.

\section{Notations, terminologies and preliminaries}
\subsection{General notations and elementary properties}
\begin{enumerate}
\item[(1)] If $X$ is a finite set, we let $\nb{X}$ denote the
cardinality of $X$. If in addition $X \subset \mathbb{Z}$ and $X
\neq \emptyset$, we let $\diam(X)$ denote the usual diameter of
$X$ (that is $\diam(X) := \max_{x , y \in X} \nb{x - y}$) and we
let
$$\delta(X) := \gcd\{x - y ~|~ x , y \in X\}$$
(with the convention $\delta(X) = 1$ if $\nb{X} = 1$).
 \item[(2)] If $B$ and $C$ are two sets of integers, the
notation $B \sim C$ means that the symmetric difference $B \Delta
C$ $(= (B \setminus C) \cup (C \setminus B))$ is finite; namely
$B$ and $C$ differ just by a finite number of elements.
 \item[(3)] If $A_1 , A_2 , \dots , A_n$ $(n \geq 1)$ are nonempty subsets of an abelian
group, we write
$$\sum_{i = 1}^{n} A_i := \{a_1 + a_2 + \dots + a_n ~|~ a_i \in A_i\} .$$
If $A_1 = A_2 = \dots = A_n \neq \mathbb{Z}$, it is convenient to
write the previous set as $n A_1$; while $n \mathbb{Z}$ stands for
the set of the integer multiples of $n$.
 \item[(4)] If $U = {(u_i)}_{i \in \mathbb{N}}$ is a nondecreasing and
 non-stationary sequence of integers, we write, for all $m \in
 \mathbb{N}$, $U(m)$ the number of terms of $U$ not exceeding
 $m$.\\
 (Stress that if $U$ is increasing, then it is just considered as a subset of $\mathbb{Z}$ having
 a finite intersection with $\mathbb{Z}^-$).\\
 $\bullet$ We call ``the lower asymptotic density'' of $U$ the quantity
 defined by
 $$\d(U) := \liminf_{m \rightarrow + \infty} \frac{U(m)}{m} \in [0 , + \infty] .$$
If $U$ is increasing (so it is a subset of $\mathbb{Z}$ having a
finite intersection with $\mathbb{Z}^-$), we clearly have $\d(U)
\leq 1$.
 \item[(5)] If $U_1 , U_2 , \dots , U_n$ $(n \geq 1)$ are
 nondecreasing and non-stationary sequences of integers, indexed
 by $\mathbb{N}$, the notation $U_1 \vee U_2 \vee \dots \vee
 U_n$ (or $\vee_{i = 1}^{n} U_i$) represents the aggregate of the
 elements of $U_1 , \dots , U_n$; each element being counted
 according to its multiplicity.\\
$\bullet$ It's clear that for all $m \in \mathbb{N}$, we have:
$(U_1 \vee \dots \vee U_n)(m) = \sum_{i = 1}^{n} U_i(m)$. So, it
follows that:
$$
\d(U_1 \vee \dots \vee U_n) \geq \sum_{i = 1}^{n} \d(U_i) .
$$
$\bullet$ Further, if $U_1 , \dots , U_n$ are increasing (so they
are simply sets), we clearly have:
$$
\d(U_1 \vee \dots \vee U_n) \geq \d(U_1 \cup \dots \cup U_n) .
$$
\item[(6)] It is easy to check that if $U$ is a nondecreasing and
non-stationary sequence of integers (indexed by $\mathbb{N}$) and
$t \in \mathbb{Z}$, then we have:
$$
(U + t)(m) = U(m) + O(1) .
$$
\item[(7)] If $B$ is a nonempty set of integers and $g$ is a
positive integer, we denote $\frac{B}{g \mathbb{Z}}$
 the image of $B$ under the canonical surjection $\mathbb{Z} \rightarrow \frac{\mathbb{Z}}{g
 \mathbb{Z}}$. We also denote $B^{(g)}$ the set of all natural
 numbers which are congruent modulo $g$ to some element of $B$; in other words:
 $$B^{(g)} := (B + g \mathbb{Z}) \cap \mathbb{N} .$$
 $\bullet$ We can easily check that if $B$ and $C$ are two
 nonempty sets of integers and $g$ is a positive integer, then we
 have:
 $$
(B + C)^{(g)} \sim B^{(g)} + C .
 $$
 In particular, if we have $B \sim B^{(g)}$ then we also have $B + C \sim (B +
 C)^{(g)}$.
\end{enumerate}
\subsection{The theorems of Kneser (see \cite{hr}, Chap 1)}
\begin{thm}[The first theorem of Kneser]~\\
Let $A_1 , A_2 , \dots , A_n$ $(n \geq 1)$ be nonempty sets of
integers having each one a finite intersection with
$\mathbb{Z}^-$. Then either
\begin{equation}\label{eq5}
\d\left(\sum_{i = 1}^{n} A_i\right) \geq \d\left(\bigvee_{i =
1}^{n} A_i\right) \tag{I}
\end{equation}
or there exists a positive integer $g$ such that
\begin{equation}\label{eq6}
\sum_{i = 1}^{n} A_i \sim \left(\sum_{i = 1}^{n} A_i\right)^{(g)}
. \tag{II}
\end{equation}
\end{thm}
{\bf Remarks:}\\
$\bullet$ We call (\ref{eq5}) ``the first alternative of the first
theorem of Kneser'' and we call (\ref{eq6}) ``the second
alternative of the first theorem of Kneser''.\\
$\bullet$ The relation (\ref{eq6}) implies in particular that the
set $\sum_{i = 1}^{n} A_i$ is (starting from some element) a
finite union of arithmetic progressions with common difference
$g$.
\begin{thm}[The second theorem of Kneser]\label{t2}~\\
Let $G$ be a finite abelian group and $B$ and $C$ be two nonempty
subsets of $G$. Then, there exists a subgroup $H$ of $G$ such that
$$B + C = B + C + H$$
and
$$\nb{B + C} \geq \nb{B + H} + \nb{C + H} - \nb{H} .$$
\end{thm}
In the applications, we use the second theorem of Kneser in the
form given by the corollary below. We first need to define the
so-called ``a subset not degenerate of an abelian group'' and then
to give a simple property related to this one.~\vspace{1mm}\\
{\bf Definitions:}\\
$\bullet$ If $G$ is an abelian group and $B$ is a subset of $G$,
we say that ``$B$ is not degenerate in $G$'' if we have
$\stab_G(B) = \{0\}$ (where $\stab_G(B)$ denotes the stabilizer of $B$ in $G$).\\
$\bullet$ If $B$ is a set of integers and $g$ is a positive
integer, we say that ``$B$ is not degenerate modulo $g$'' if
$\frac{B}{g \mathbb{Z}}$ is not degenerate in $\frac{\mathbb{Z}}{g
\mathbb{Z}}$.~\vspace{1mm}\\
\begin{prop}\label{p2}
Let $G$ be an abelian group and $B$ and $C$ be two nonempty
subsets of $G$ such that $(B + C)$ is not degenerate in $G$. Then
also $B$ and $C$ are not degenerate in $G$.
\end{prop}
{\bf Proof.} This is an immediate consequence of the fact that:\\
$\stab_G(B) + \stab_G(C) \subset \stab_G(B +
C)$.\penalty-20\null\hfill$\blacksquare$\par\medbreak
\begin{coll}\label{c1}
Let $G$ be a finite abelian group and $B_1 , \dots , B_n$ $(n \geq
1)$ be nonempty subsets of $G$ such that $(B_1 + \dots + B_n)$ is
not degenerate in $G$. Then we have
$$\nb{B_1 + \dots + B_n} \geq \nb{B_1} + \dots + \nb{B_n} - n + 1 .$$
\end{coll}
{\bf Proof.} It suffices to show the corollary for $n = 2$. The
general case follows by a simple induction on $n$ and by using
Proposition \ref{p2}. Suppose $n= 2$. Theorem \ref{t2} gives a
subgroup $H$ of $G$ satisfying the two relations $B_1 + B_2 = B_1
+ B_2 + H$ and $\nb{B_1 + B_2} \geq \nb{B_1 + H} + \nb{B_2 + H} -
\nb{H}$. The first one implies $H \subset \stab_G(B_1 + B_2) =
\{0\}$, so $H = \{0\}$. By replacing this into the second one, we
conclude to $\nb{B_1 + B_2} \geq \nb{B_1} + \nb{B_2} - 1$ as
required.\penalty-20\null\hfill$\blacksquare$\par\medbreak
 The following proposition (which is an easy exercise) makes the
 connection between the first and the second theorem of Kneser:
 \begin{prop}\label{p1}
Let B be a nonempty set of integers and $g$ be a positive integer.
The two following assertions are equivalent:
\begin{enumerate}
\item[{\rm (i)}] $B$ is not degenerate modulo $g$
 \item[{\rm (ii)}] There is no positive integer $m < g$ such that $B^{(m)} =
 B^{(g)}$.
\end{enumerate}
 \end{prop}

 Now, let us explain how we use the theorems of Kneser in this paper.
We first get sets $A_i = h_i (A \setminus X)$, $i = 0 , \dots ,
 n$ such that $\cup_{i = 1}^{n} (A_i + \tau_i) \sim \mathbb{N}$ and $\d(A_0) > 0$
 (where $n$ is a natural number depending on $A$ and $X$, the $h_i$'s are positive integers
 depending only on $h$ and such that $h_0 \leq n$ and the $\tau_i$'s are integers). We
thus have $\d(\vee_{i = 0}^{n} A_i) > 1$, implying that
 the first alternative of the first theorem of Kneser cannot hold.
 Consequently we are in the second alternative of the first
 theorem of Kneser, namely there exists a positive integer $g$
 such that $\sum_{i = 0}^{n} A_i \sim (\sum_{i = 0}^{n}
 A_i)^{(g)}$. By choosing $g$ minimal to have this property, we deduce from Proposition \ref{p1} that the set
 $\sum_{i = 0}^{n} A_i$ is not degenerate modulo $g$; in other
 words the set $\sum_{i = 0}^{n} \frac{A_i}{g \mathbb{Z}}$ is not
 degenerate in the group $\frac{\mathbb{Z}}{g \mathbb{Z}}$. It follows from Proposition \ref{p2} that also
 $\sum_{i = 1}^{n} \frac{A_i}{g \mathbb{Z}}$ is not degenerate in $\frac{\mathbb{Z}}{g \mathbb{Z}}$. Then
 by applying Corollary \ref{c1} for $G = \frac{\mathbb{Z}}{g
 \mathbb{Z}}$ and $B_i = \frac{A_i}{g \mathbb{Z}}$ $(i = 1 , \dots ,
 n)$, we deduce that $\nb{\frac{\sum_{i = 1}^{n} A_i}{g \mathbb{Z}}} \geq \sum_{i = 1}^{n}
 \nb{\frac{A_i}{g \mathbb{Z}}} - n + 1 \geq g - n + 1$ (since $\cup_{i = 1}^{n} (A_i + \tau_i) \sim
 \mathbb{N}$); so $\nb{\frac{(h_1 + \dots + h_n) (A \setminus X)}{g \mathbb{Z}}} \geq g - n +
 1$. Next, from the nature of the sequence ${(\nb{\frac{r (A \setminus X)}{g
 \mathbb{Z}}})}_{r \in \mathbb{N}}$ (pointed out in Lemma \ref{l1} of the next
 section) and the hypothesis that $A \setminus X$ is a basis, we derive that
 $\nb{\frac{(h_1 + \dots + h_n + n) (A \setminus X)}{g \mathbb{Z}}} =
 g$; hence $\frac{(h_1 + \dots + h_n + n) (A \setminus X)}{g \mathbb{Z}} = \frac{\mathbb{Z}}{g
 \mathbb{Z}}$. We thus have $((h_1 + \dots + h_n + n)(A \setminus X))^{(g)} \sim
 \mathbb{N}$. But since on the other hand we have (in view of the elementary properties of §2.1):
 $((h_1 + \dots + h_n + n)(A \setminus X))^{(g)}
 = ((A_0 + \dots + A_n) + (n - h_0)(A \setminus X))^{(g)} \sim (A_0 + \dots + A_n)^{(g)} + (n - h_0)(A \setminus X)
 \sim A_0 + \dots + A_n + (n - h_0)(A \setminus X) = (h_1 + \dots + h_n + n) (A \setminus
 X)$, it finally follows that $(h_1 + \dots + h_n + n) (A \setminus X) \sim
 \mathbb{N}$, that is $\G(A \setminus X) \leq h_1 + \dots + h_n +
 n$.

 In the work of Nash \cite{n}, the parameter $n$ depends on $h$
 and $\nb{X}$. Actually, its dependence in $\nb{X}$ stems from
 the upper bounds of the cardinalities of the sets
 $\ell X$ $(\ell = 0 , \dots , h)$. In \cite{n}, the upper bound used for each $\nb{\ell
 X}$ is $\binom{\nb{X} + \ell - 1}{\ell}$, which is a polynomial in $\ell$
with degree $(\nb{X} - 1)$ and
 then leads to bound from above $\G(A \setminus X)$ by a polynomial in $h$ with degree $(\nb{X} + 1)$.
 However, that estimate of $\nb{\ell X}$ is very large for many sets $X$; for example if
 $X$ is an arithmetic progression, we simply have $\nb{\ell X} = \ell \nb{X} - \ell +
 1$ which is linear in $\ell$ and (as we will see it later) allows to estimate $\G(A \setminus
 X)$ by a polynomial with degree $3$ in $h$. In order to obtain such an estimate for $\G(A \setminus X)$
 in the general case, our idea (see Lemmas \ref{l2} and \ref{l3}) consists to
 replace $\nb{X}$ by another parameter in
 $X$ (resp. $X$ and $A$) for which the cardinality of each of the sets $\ell X$
 (resp. other more complex sets) is bounded above by a
 linear function in $\ell$ (resp. simple function in $h$). The upper bounds obtained in this way for
 $\G(A \setminus X)$ are simply polynomials in $h$ with degrees $3$ or $2$ and with coefficients
 linear in the considered parameters (see Theorems \ref{t3} and
\ref{t4}). On the other hand,
 it must be noted that upper bounds for $\G(A \setminus X)$ which
 are polynomials with degrees $3$ or $2$ in $h$ can be directly
 derived from the theorem of Nash, but in this way we lose the linearity in the considered
 parameter (see Theorem \ref{t5} and Remark \ref{r1}).
\section{Lemmas}~

The two first lemmas which follow constitute the main differences
with Nash' work \cite{n} about the use of Kneser's theorems. While
the third one gives the nature (in terms of monotony) of some
sequences (related to a given finite abelian group) which also
plays a vital part in the proof of our results.
\begin{lemma}\label{l2}
Let $X$ be a nonempty finite set of integers. Then we have:
$$\nb{X} \leq \frac{\diam(X)}{\delta(X)} + 1 .$$
In addition, this inequality becomes an equality if and only if
$X$ is an arithmetic progression.
\end{lemma}
{\bf Proof.} The lemma is obvious if $\nb{X} = 1$. Assume for the
following that $\nb{X} \geq 2$ and write $X = \{x_0 , x_1 , \dots
, x_n\}$ $(n \geq 1)$, with $x_0 < x_1 < \dots < x_n$. Since the
positive integers $x_i - x_{i - 1}$ $(i = 1 , \dots , n)$ are
clearly multiples of $\delta(X)$ then we have $x_i - x_{i - 1}
\geq \delta(X)$ $(\forall i = 1 , \dots , n)$. It follows that
$x_n - x_0 = \sum_{i = 1}^{n} (x_i - x_{i - 1}) \geq n \delta(X)$,
which gives $n \leq \frac{x_n - x_0}{\delta(X)} =
\frac{\diam(X)}{\delta(X)}$. Hence $\nb{X} = n + 1 \leq
\frac{\diam(X)}{\delta(X)} + 1$ as required.\\
Further, the above proof shows well that the inequality of the
lemma is reached if and only if we have $x_i - x_{i - 1} =
\delta(X)$ $(\forall i = 1 , \dots , n)$ which simply means that
$X$ is an arithmetic progression. The proof is
complete.\penalty-20\null\hfill$\blacksquare$\par\medbreak

\begin{lemma}\label{l3}
Let $X$ be a finite nonempty set of integers and $B$ be an
infinite set of integers having a finite intersection with
$\mathbb{Z}^-$. Define:
$$\eta := \min_{\begin{array}{c}\scriptstyle{b , b' \in B , b \neq b'} \\ \scriptstyle{\nb{b - b'} \geq
\diam(X)}\end{array}}\nb{b - b'} .$$ Then, for all $u , v \in
\mathbb{N}$, $g \in \mathbb{N}^*$, we have:
$$(u B + v X)(m) \leq \eta . ((u + v) B)(m) + O(1)$$
and
$$\nb{\frac{u B + v X}{g \mathbb{Z}}} \leq \eta \nb{\frac{(u + v) B}{g \mathbb{Z}}} .$$
\end{lemma}
{\bf Proof.} Since we have for all $\tau \in \mathbb{Z}$: $(u B +
v X + \tau)(m) = (u B + v X)(m) + O(1)$ (according to the part (6)
of §2.1) and $\nb{\frac{u B + v X + \tau}{g \mathbb{Z}}} =
\nb{\frac{u B + v X}{g \mathbb{Z}}}$ (obviously), then there is no
loss of generality in translating $B$ and $X$ by integers. By
translating, if necessary, $X$, assume that $0$ is its smaller
element and write $X = \{x_0 , x_1 , \dots , x_n\}$ $(n \in
\mathbb{N})$, with $0 = x_0 < x_1 < \dots < x_n$. Next, let $b_0 ,
b \in B$ such that $b - b_0 = \eta$. By translating, if necessary,
$B$, assume $b_0 = 0$. Then we have
$$b = \eta \geq \diam(X) = x_n .$$
In this situation, we claim that we have
\begin{equation}\label{eq9}
(u B + v X) \subset \bigcup_{0 \leq \tau < \eta}\left((u + v) B +
\tau\right)
\end{equation}
which clearly implies the two inequalities of the lemma. So, it
just remains to show (\ref{eq9}). Let $N \in (u B + v X)$ and show
that there exists a non-negative integer $\tau < \eta$ such that
$N \in (u + v) B + \tau$. Since $0 = b_0 = x_0 \in B \cap X$, the
fact that $N \in (u B + v X)$ means that $N$ can be written in the
form
\begin{equation}\label{eq10}
N = u_1 b_1 + \dots + u_k b_k + v_1 x_1 + \dots + v_n x_n ,
\end{equation}
with $k , u_1 , \dots , u_k , v_1 , \dots , v_n \in \mathbb{N}$,
$b_1 , \dots , b_k \in B$, $u_1 + \dots + u_k \leq u$ and $v_1 +
\dots + v_n \leq v$.\\
Now, since $x_1 < x_2 < \dots < x_n \leq \eta$, then we have $v_1
x_1 + \dots + v_n x_n \leq (v_1 + \dots + v_n) \eta \leq v \eta$,
which implies that the euclidean division of the non-negative
integer $(v_1 x_1 + \dots + v_n x_n)$ by $\eta$ yields:
\begin{equation}\label{eq11}
v_1 x_1 + \dots + v_n x_n = t \eta + \tau ,
\end{equation}
with $t , \tau \in \mathbb{N}$, $t \leq v$ and $0 \leq \tau <
\eta$. By reporting (\ref{eq11}) into (\ref{eq10}), we finally
obtain
\begin{equation}\label{eq12}
N = u_1 b_1 + \dots + u_k b_k + t \eta + \tau .
\end{equation}
Since $0 = b_0 \in B$, $b_1 , \dots , b_k , \eta \in B$ (recall
that $\eta = b$) and $u_1 + \dots + u_k + t \leq u + v$, then the
relation (\ref{eq12}) is well a writing of $N$ as a sum of $(u +
v)$ elements of $B$ and $\tau$; in other words $N \in (u + v) B +
\tau$, giving the desired conclusion. The proof is
complete.\penalty-20\null\hfill$\blacksquare$\par\medbreak

\begin{lemma}\label{l1}
Let $G$ be a finite abelian group and $B$ be a nonempty subset of
$G$. For all $r \in \mathbb{N}$, set $u_r := |r B|$. Then, there
exists $r_0 \in \mathbb{N}$ such that:
$$u_0 < u_1 < \dots < u_{r_0}$$
and
$$~~~~~~~~ u_r = u_{r_0} ~~~~~~~~ (\forall r \geq r_0) .$$
\end{lemma}
{\bf Proof.} Firstly, since $G$ is finite, the sequence
${(u_r)}_r$ is bounded above by $\nb{G}$. Secondly, we claim that
${(u_r)}_r$ is nondecreasing. Indeed, by fixing $b \in B$, we have
for all $r \in \mathbb{N}$: $(r + 1) B \supset r B + b$, hence
$u_{r + 1} = \nb{(r + 1) B} \geq \nb{r B + b} = \nb{r B} = u_r$.
It follows from these two facts that there exists $r_0 \in
\mathbb{N}$ such that $u_{r_0} = u_{r_0 + 1}$. By taking $r_0$
minimal to have this property, we have:
$$u_0 < u_1 < \dots < u_{r_0} = u_{r_0 + 1} .$$
To conclude the proof of the lemma, it remains to show that
\begin{equation}\label{eq7}
u_r = u_{r_0} ~~~~~~ (\forall r \geq r_0) .
\end{equation}
If $b \in B$ is fixed, we claim that for all $r \geq r_0$, we
have:
\begin{equation}\label{eq8}
r B = r_0 B + (r - r_0) b
\end{equation}
which clearly implies (\ref{eq7}). So, it remains to show
(\ref{eq8}). To do this, we argue by induction on $r \geq r_0$.
For $r = r_0$, the relation (\ref{eq8}) is obvious. Next, since
$(r_0 + 1) B \supset r_0 B + b$ and $\nb{(r_0 + 1) B} = u_{r_0 +
1} = u_{r_0} = \nb{r_0 B} = \nb{r_0 B + b}$, then we certainly
have $(r_0 + 1) B = r_0 B + b$, showing that (\ref{eq8}) also
holds for $r = r_0 + 1$. Now, let $r \geq r_0$, assume that
(\ref{eq8}) holds for $r$ and show that it also holds for $(r +
1)$. We have:
\begin{eqnarray*}
(r + 1) B & = & (r_0 + 1) B + (r - r_0) B \\
& \!\!\!\!\!\!\!\!\!\!= & \!\!\!\!\!\!\!(r_0 B + b) + (r
- r_0) B ~~~~ \text{(since (\ref{eq8}) holds for $(r_0 + 1)$)} \\
& \!\!\!\!\!\!\!\!\!\!= & \!\!\!\!\!\!\!r B + b \\
& \!\!\!\!\!\!\!\!\!\!= & \!\!\!\!\!\!\!\left(r_0 B + (r
- r_0) b\right) + b ~~~~ \text{(from the induction hypothesis)} \\
& \!\!\!\!\!\!\!\!\!\!= & \!\!\!\!\!\!\!r_0 B + (r + 1 - r_0) b .
\end{eqnarray*}
Hence (\ref{eq8}) also holds for $(r + 1)$. This finishes this
induction and completes the
proof.\penalty-20\null\hfill$\blacksquare$\par\medbreak

\section{Main Results}~

Throughout this section, we fix an additive basis $A$ and a finite
nonempty subset $X$ of $A$ such that $A \setminus X$ is still a
basis. We put $h := \G(A)$ and we define
$$
d := \frac{\diam(X)}{\delta(X)} ~~,~~ \eta :=
\min_{\begin{array}{c}\scriptstyle{a , b \in A \setminus X , a
\neq b}
\\ \scriptstyle{\nb{a - b} \geq \diam(X)}\end{array}} \!\!\!\!\!\!\nb{a -
b} ~~\text{and}~~ \mu := \min_{y \in A \setminus X} \diam(X \cup
\{y\}) .$$
\begin{thm}\label{t3}
We have $\displaystyle \G(A \setminus X) \leq \frac{h (h + 3)}{2}
+ d \frac{h (h - 1) (h + 4)}{6}$.
\end{thm}
{\bf Proof.} Put $B := A \setminus X$, so $A = B \cup X$. Then,
the fact that $A$ is a basis of order $h$ amounts to:
\begin{equation}\label{eq13}
h B \cup ((h - 1) B + X) \cup ((h - 2) B + 2 X) \cup \dots \cup (B
+ (h - 1) X) \sim \mathbb{N} .
\end{equation}
(Remark that $h X$ is finite).\\
Now, since the set of the left-hand side of (\ref{eq13}) is
clearly contained in a finite union of translates of $h B$, then
by denoting $N$ a number of translates of $h B$ which are
sufficient to cover it, we have (according to the part (6) of
§2.1):
$$\left(h B \cup ((h - 1) B + X) \cup \dots \cup (B + (h - 1) X)\right)(m) \leq N . (h B)(m) + O(1) .$$
It follows that:
\begin{eqnarray*}
\liminf_{m \rightarrow + \infty} \frac{(h B)(m)}{m} & ~ & \\ &
\!\!\!\!\!\!\!\!\!\!\!\!\!\!\!\!\!\!\!\!\!\!\!\!\!\!\!\!\!\!\!\!\!\!\!\!\!\!\!\!\!\!\!\!\!\!\!\!\!\!\!\!\geq
&
\!\!\!\!\!\!\!\!\!\!\!\!\!\!\!\!\!\!\!\!\!\!\!\!\!\!\!\!\frac{1}{N}
\liminf_{m \rightarrow + \infty} \frac{1}{m}\left(h B \cup ((h -
1) B + X) \cup \dots \cup (B + (h
- 1) X)\right)(m) \\
&
\!\!\!\!\!\!\!\!\!\!\!\!\!\!\!\!\!\!\!\!\!\!\!\!\!\!\!\!\!\!\!\!\!\!\!\!\!\!\!\!\!\!\!\!\!\!\!\!\!\!\!\!=
&
\!\!\!\!\!\!\!\!\!\!\!\!\!\!\!\!\!\!\!\!\!\!\!\!\!\!\!\!\frac{1}{N}
~~~~~~ \text{(according to (\ref{eq13}))} .
\end{eqnarray*}
Thus
\begin{equation}\label{eq14}
\d(h B) \geq \frac{1}{N} > 0 .
\end{equation}
Now, according to (\ref{eq13}), (\ref{eq14}) and the part (5) of
§2.1, we have:
\begin{eqnarray*}
& ~ & \hspace{-1cm}\d\left(h B \vee h B \vee ((h - 1) B + X) \vee
((h - 2) B + 2 X)
\vee \dots \vee (B + (h - 1) X)\right) \\
& \geq & \d(h B) + \d\left(h B \vee ((h - 1) B + X) \vee \dots
\vee (B + (h - 1) X)\right) \\
& \geq & \d(h B) + \d\left(h B \cup ((h - 1) B + X) \cup \dots
\cup (B + (h - 1) X)\right) \\
& = & \d(h B) + 1 \\
& > & 1 .
\end{eqnarray*}
So, we have
\begin{equation}\label{eq15}
\begin{split}
\liminf_{m \rightarrow + \infty} \frac{1}{m} \{& (h
B)(m) + (h B)(m) + ((h - 1) B + X)(m) \\
&+ ((h - 2) B + 2 X)(m) + \dots + (B + (h - 1) X)(m)\}
> 1 .
\end{split}
\end{equation}
Next, according to the part (6) of §2.1 and to Lemma \ref{l2},
each of the quantities $((h - \ell) B + \ell X)(m)$ $(\ell = 1 ,
\dots , h - 1)$ is bounded above as follows
\begin{equation}\label{eq16}
\begin{split}
((h - \ell) B + \ell X)(m) &\leq \nb{\ell X} . ((h - \ell) B)(m) +
O(1) \\
&\leq \left(\frac{\diam(\ell X)}{\delta(\ell X)} + 1\right) . ((h
-
\ell) B)(m) + O(1) \\
&= (\ell d + 1) . ((h - \ell) B)(m) + O(1)
\end{split}
\end{equation}
(since $\diam(\ell X) = \ell \diam(X)$ and $\delta(\ell X) =
\delta(X)$). \\ Then, by reporting these into (\ref{eq15}), we
obtain:
\begin{equation*}
\begin{split}
\liminf_{m \rightarrow + \infty} \frac{1}{m} \{ &(h B)(m) + (h
B)(m) + (d + 1).((h - 1) B)(m) \\ &+ (2 d + 1).((h - 2) B)(m) +
\dots + ((h - 1)d + 1).B(m)\} > 1 ,
\end{split}
\end{equation*}
which amounts to
\begin{equation}\label{eq17}
\d\left(h B \vee \bigvee_{\ell = 0}^{h -
1}\left(\bigvee_{\text{$(\ell d + 1)$ times}} (h - \ell)
B\right)\right) > 1 .
\end{equation}
This last relation shows well that the first alternative of the
first theorem of Kneser (applied to the set $h B$ with $(\ell d +
1)$ copies of each of the sets $(h - \ell) B$, $\ell = 0 , 1 ,
\dots , h - 1$) cannot hold. We are thus in the second alternative
of the first theorem of Kneser; that is there exists a positive
integer $g$ such that
\begin{equation}\label{eq18}
\left(h + \sum_{\ell = 0}^{h - 1} (\ell d + 1)(h - \ell)\right) B
\sim \left(\left(h + \sum_{\ell = 0}^{h - 1} (\ell d + 1)(h -
\ell)\right) B\right)^{(g)} .
\end{equation}
Let's take $g$ minimal in (\ref{eq18}). This implies from
Proposition \ref{p1} that the set $(h + \sum_{\ell = 0}^{h - 1}
(\ell d + 1)(h - \ell)) B$ is not degenerate modulo $g$; in other
words, the set $(h + \sum_{\ell = 0}^{h - 1} (\ell d + 1) (h -
\ell)) \frac{B}{g \mathbb{Z}}$ is not degenerate in
$\frac{\mathbb{Z}}{g \mathbb{Z}}$. It follows from Proposition
\ref{p2} that also the set $(\sum_{\ell = 0}^{h - 1} (\ell d + 1)
(h - \ell)) \frac{B}{g \mathbb{Z}}$ is not degenerate in
$\frac{\mathbb{Z}}{g \mathbb{Z}}$. Then, from Corollary \ref{c1},
we have
\begin{eqnarray}
\!\!\!\!\!\!\!\!\!\!\!\!\!\!\!\!\!\!\!\!\!\!\!\!\!\!\!\!\!\!\!\!\!\!\!\!\!\!\!\!\nb{\left(\sum_{\ell
= 0}^{h - 1} (\ell d + 1)(h - \ell)\right) \frac{B}{g \mathbb{Z}}}
& = & \nb{\sum_{\ell = 0}^{h - 1} \sum_{\text{$(\ell d + 1)$
times}} \frac{(h - \ell) B}{g \mathbb{Z}}} \notag \\
&
\!\!\!\!\!\!\!\!\!\!\!\!\!\!\!\!\!\!\!\!\!\!\!\!\!\!\!\!\!\!\!\!\!\!\!\!\!\!\!\!\!\!\!\!\!\!\!\!\!\!
\!\!\!\!\!\!\!\!\!\!\!\!\!\!\!\!\!\!\!\!\!\!\!\!\!\!\!\!\!\!\geq &
\!\!\!\!\!\!\!\!\!\!\!\!\!\!\!\!\!\!\!\!\!\!\!\!\!\!\!\!\!\!\!\!\!\!\!\!\!\!\!\!\sum_{\ell
= 0}^{h - 1} (\ell d + 1) \nb{\frac{(h - \ell) B}{g \mathbb{Z}}} -
\sum_{\ell = 0}^{h - 1} (\ell d + 1) + 1 . \label{eq19}
\end{eqnarray}
Now, let's bound from below the sum $\sum_{\ell = 0}^{h - 1} (\ell
d + 1) \nb{\frac{(h - \ell) B}{g \mathbb{Z}}}$. We have for all
$\ell \in \{0 , 1 , \dots , h - 1\}$:
\begin{eqnarray*}
(\ell d + 1) \nb{\frac{(h - \ell) B}{g \mathbb{Z}}} & = &
\left(\frac{\diam(\ell X)}{\delta(\ell X)} + 1\right) \nb{\frac{(h
- \ell) B}{g \mathbb{Z}}} \\
& \geq & \nb{\ell X} . \nb{\frac{(h - \ell) B}{g \mathbb{Z}}}
~~~~~~~~\text{(according to Lemma \ref{l2})} \\
& \geq & \nb{\frac{\ell X}{g \mathbb{Z}}} . \nb{\frac{(h - \ell)
B}{g \mathbb{Z}}} \\
& \geq & \nb{\frac{(h - \ell) B + \ell X}{g \mathbb{Z}}} ;
\end{eqnarray*}
hence
\begin{eqnarray*}
\sum_{\ell = 0}^{h - 1} (\ell d + 1) \nb{\frac{(h - \ell) B}{g
\mathbb{Z}}} & \geq & \sum_{\ell = 0}^{h - 1} \nb{\frac{(h - \ell)
B + \ell X}{g \mathbb{Z}}} \\
& \geq & \nb{\frac{h B \cup ((h - 1) B + X) \cup \dots \cup (B +
(h - 1) X)}{g \mathbb{Z}}} \\
& = & g ~~~~~~~~~~~~~~~~~~~~~~\text{(according to (\ref{eq13}))} .
\end{eqnarray*}
By reporting this into (\ref{eq19}), we have
\begin{equation}\label{eq20}
\nb{\left(\sum_{\ell = 0}^{h - 1} (\ell d + 1)(h - \ell)\right)
\frac{B}{g \mathbb{Z}}} \geq g - \sum_{\ell = 0}^{h - 1} (\ell d +
1) + 1 .
\end{equation}
Now, from Lemma \ref{l1}, we know that the sequence of natural
numbers ${(\nb{r \frac{B}{g \mathbb{Z}}})}_{r \in \mathbb{N}}$
increases until reaching its maximal value which it then continues
to take indefinitely. In addition, because $\G(B) B \sim
\mathbb{N}$, we have $\nb{\G(B) \frac{B}{g \mathbb{Z}}} =
\nb{\frac{\mathbb{Z}}{g \mathbb{Z}}} = g$, showing that $g$ is the
maximal value of the same sequence. On the other hand, if we
assume that the finite sequence \\${(\nb{r \frac{B}{g
\mathbb{Z}}})}_{\sum_{\ell = 0}^{h - 1} (\ell d + 1) (h - \ell)
\leq r \leq \sum_{\ell = 0}^{h - 1} (\ell d + 1) (h - \ell + 1)}$
is increasing, we would have (according to (\ref{eq20})):
$$\nb{\left(\sum_{\ell = 0}^{h - 1} (\ell d + 1) (h - \ell + 1)\right) \frac{B}{g \mathbb{Z}}} \geq g + 1$$
which is impossible. Consequently, the sequence ${(\nb{r
\frac{B}{g \mathbb{Z}}})}_{r \in \mathbb{N}}$ becomes constant
(equal to $g$) before its term of order $r = \sum_{\ell = 0}^{h -
1} (\ell d + 1) (h - \ell + 1)$. In particular, we have
$$\nb{\left(\sum_{\ell = 0}^{h - 1} (\ell d + 1)(h - \ell + 1)\right) \frac{B}{g \mathbb{Z}}} = g$$
and then
$$\left(\sum_{\ell = 0}^{h - 1} (\ell d + 1)(h - \ell + 1)\right)\frac{B}{g \mathbb{Z}} =
\frac{\mathbb{Z}}{g \mathbb{Z}} ,$$ implying that
\begin{equation}\label{eq21}
\left(\left(\sum_{\ell = 0}^{h - 1} (\ell d + 1)(h - \ell +
1)\right) B\right)^{(g)} = \mathbb{N} .
\end{equation}
But on the other hand, since $\sum_{\ell = 0}^{h - 1} (\ell d + 1)
(h - \ell + 1) \geq h + \sum_{\ell = 0}^{h - 1} (\ell d + 1) (h -
\ell)$, we have (according to the relation (\ref{eq18}) and the
property of the part (7) of §2.1):
\begin{equation}\label{eq22}
\left(\sum_{\ell = 0}^{h - 1} (\ell d + 1)(h - \ell + 1)\right) B
\sim \left(\left(\sum_{\ell = 0}^{h - 1} (\ell d + 1)(h - \ell +
1)\right) B\right)^{(g)} .
\end{equation}
By comparing (\ref{eq21}) and (\ref{eq22}), we finally deduce that
$$\left(\sum_{\ell = 0}^{h - 1} (\ell d + 1)(h - \ell + 1)\right) B \sim \mathbb{N} ,$$
which gives
$$\G(B) \leq \sum_{\ell = 0}^{h - 1} (\ell d + 1) (h - \ell + 1) = \frac{h(h + 3)}{2} + d \frac{h(h - 1)(h + 4)}{6}$$
(since $\sum_{\ell = 0}^{h - 1} \ell = \frac{h(h - 1)}{2}$ and
$\sum_{\ell = 0}^{h - 1} \ell^2 = \frac{h(h - 1)(2 h - 1)}{6}$).\\
The theorem is
proved.\penalty-20\null\hfill$\blacksquare$\par\medbreak
\begin{coll}\label{c2}
If in addition $X$ is an arithmetic progression, then we have:
$$\G(A \setminus X) \leq \frac{h (h + 3)}{2} + (\nb{X} - 1) \frac{h (h - 1) (h + 4)}{6} .$$
\end{coll}
{\bf Proof.} By Lemma \ref{l2}, we have $\nb{X} =
\frac{\diam(X)}{\delta(X)} + 1 = d + 1$, hence $d = \nb{X} - 1$.
The corollary then follows at once from Theorem
\ref{t3}.\penalty-20\null\hfill$\blacksquare$\par\medbreak
\begin{thm}\label{t4}
We have $\displaystyle \G(A \setminus X) \leq \eta (h^2 - 1) + h +
1$.
\end{thm}
{\bf Proof.} We proceed as in the proof of Theorem \ref{t3} with
some differences; so we only detail these differences. Putting $B
:= A \setminus X$, we repeat the proof of Theorem \ref{t3} until
the relation (\ref{eq15}). After that, using Lemma \ref{l3}, we
bound from above each of the quantities $((h - \ell) B + \ell
X)(m)$ $(\ell = 1 , \dots , h - 1)$ by
\begin{equation}\label{eq23}
((h - \ell) B + \ell X)(m) \leq \eta . (h B)(m) + O(1) .
\tag{\ref{eq16}$'$}
\end{equation}
Then, by reporting these into (\ref{eq15}), we obtain
\begin{equation}\label{eq24}
\d\left(\bigvee_{\text{$(\eta (h - 1) + 2)$ times}} (h B)\right) >
1 , \tag{\ref{eq17}$'$}
\end{equation}
which shows well that the first alternative of the first theorem
of Kneser (applied to $(\eta (h - 1) + 2)$ copies of the set $h
B$) cannot hold. Consequently, we are in the second alternative of
the first theorem of Kneser, that is there exists a positive
integer $g$ such that
\begin{equation}\label{eq25}
(\eta (h - 1) + 2) h B \sim \left((\eta (h - 1) + 2) h
B\right)^{(g)} . \tag{\ref{eq18}$'$}
\end{equation}
Let's take $g$ minimal in (\ref{eq25}). Then, Propositions
\ref{p1} and \ref{p2} imply that the set $(\eta (h - 1) + 1) h
\frac{B}{g \mathbb{Z}}$ is non degenerate in $\frac{\mathbb{Z}}{g
\mathbb{Z}}$. It follows from Corollary \ref{c1} that we have
$$
\nb{(\eta (h - 1) + 1) h \frac{B}{g \mathbb{Z}}} =
\nb{\sum_{\text{$(\eta (h - 1) + 1)$ times}} \frac{h B}{g
\mathbb{Z}}} ~~~~~~~~~~~~~~~~~~~~~~~~~~~~~~~~~~~~~~~~~~
$$
\begin{equation}\label{eq26}
\geq (\eta (h - 1) + 1) \nb{\frac{h B}{g \mathbb{Z}}} - \eta (h -
1) .  \tag{\ref{eq19}$'$}
\end{equation}
Next, using the second inequality of Lemma \ref{l3}, we have
\begin{eqnarray*}
(\eta (h - 1) + 1) \nb{\frac{h B}{g \mathbb{Z}}} & = & \sum_{\ell
= 1}^{h - 1} \eta . \nb{\frac{((h - \ell) + \ell) B}{g
\mathbb{Z}}} + \nb{\frac{h B}{g \mathbb{Z}}} \\
& \geq & \sum_{\ell = 1}^{h - 1} \nb{\frac{(h - \ell) B + \ell
X}{g \mathbb{Z}}} + \nb{\frac{h B}{g \mathbb{Z}}} \\
& \geq & \nb{\bigcup_{\ell = 0}^{h - 1} \frac{((h - \ell) B + \ell
X)}{g \mathbb{Z}}} \\
& = & g ~~~~~~~~~~~~\text{(according to (\ref{eq13}))} .
\end{eqnarray*}
By reporting this into (\ref{eq26}), we have
\begin{equation}\label{eq27}
\nb{(\eta (h - 1) + 1) h \frac{B}{g \mathbb{Z}}} \geq g - \eta (h
- 1) . \tag{\ref{eq20}$'$}
\end{equation}
It follows from Lemma \ref{l1} (as we applied it in the proof of
Theorem \ref{t3}) that the sequence ${(\nb{r \frac{B}{g
\mathbb{Z}}})}_{r \in \mathbb{N}}$ is stationary in $g$ before its
term of order $r = (\eta (h - 1) + 1) (h + 1)$. In particular, we
have $\nb{(\eta (h - 1) + 1) (h + 1) \frac{B}{g \mathbb{Z}}} = g$;
hence $(\eta (h - 1) + 1) (h + 1) \frac{B}{g \mathbb{Z}} =
\frac{\mathbb{Z}}{g \mathbb{Z}}$, implying that
\begin{equation}\label{eq28}
\left((\eta (h - 1) + 1) (h + 1) B\right)^{(g)} \sim \mathbb{N} .
\tag{\ref{eq21}$'$}
\end{equation}
But on the other hand, since $\eta \geq 1$, we have $(\eta (h - 1)
+ 1) (h + 1) \geq (\eta (h - 1) + 2) h$, which implies (according
to the relation (\ref{eq25}) and the property of the part (7) of
§2.1) that
\begin{equation}\label{eq29}
(\eta (h - 1) + 1) (h + 1) B \sim \left((\eta (h - 1) + 1 ) (h +
1) B\right)^{(g)} . \tag{\ref{eq22}$'$}
\end{equation}
By comparing (\ref{eq28}) and (\ref{eq29}), we finally deduce that
$$(\eta (h - 1) + 1) (h + 1) B \sim \mathbb{N} ,$$
which gives $\G(B) \leq (\eta (h - 1) + 1)(h + 1) = \eta (h^2 - 1)
+ h + 1$, as required. The theorem is
proved.\penalty-20\null\hfill$\blacksquare$\par\medbreak
\begin{thm}\label{t5}
We have $\displaystyle \G(A \setminus X) \leq \frac{h \mu (h \mu +
3)}{2}$.
\end{thm}
{\bf Proof.} First, notice that $\mu \geq 1$ (since $X \neq
\emptyset$). Notice also that the parameters $h , \mu$ and $G(A
\setminus X)$ are still unchanged if we translate the basis $A$ by
an integer. Let $y_0 \in A \setminus X$ such that $\mu = \diam(X
\cup \{y_0\})$; so by translating if necessary $A$ by $(- y_0)$,
we can assume (without loss of generality) that $y_0 = 0$. Then
putting $X = \{x_1 , \dots , x_n\}$ $(n \geq 1)$ with $x_1 < x_2 <
\dots < x_n$, we have
\begin{equation}\label{eq30}
\mu = \diam(X \cup \{0\}) = \max\{\nb{x_1} , \nb{x_2} , \dots ,
\nb{x_n} , x_n - x_1\} .
\end{equation}
We are going to show that the set $(A \setminus X) \cup \{\pm 1\}$
is a basis of order $\leq h \mu$. The result of the theorem then
follows from the particular case `$k = 1$' of Theorem
\ref{t1} of Nash. We distinguish the three following cases:~\vspace{1mm}\\
{\bf 1\textsuperscript{st} case.} (if $x_1 \geq 0$)\\
In this case, the elements of $X$ are all non-negative. Let $N$ be
a natural number large enough that it can be written as a sum of
$h$ elements of $A$; that is
\begin{equation}\label{eq31}
N = a_1 + \dots + a_t + \alpha_1 x_1 + \dots + \alpha_n x_n ,
\end{equation}
with $t , \alpha_1 , \dots , \alpha_n \in \mathbb{N}$, $a_1 ,
\dots , a_t \in A \setminus X$ and $t + \alpha_1 + \dots +
\alpha_n = h$.\\
Next, since the non-negative integer $(\alpha_1 x_1 + \dots +
\alpha_n x_n)$ is obviously bounded above by $(\alpha_1 + \dots +
\alpha_n) \mu = (h - t) \mu \leq h \mu - t$, then it is a sum of
$(h \mu - t)$ elements of the set $\{0 , 1\}$. It follows from
(\ref{eq31}) that $N$ is a sum of $h \mu$ elements of the set $(A
\setminus X) \cup \{0 , 1\} = (A \setminus X) \cup \{1\}$. This
last fact shows well (since $N$ is an arbitrary sufficiently large
integer) that the set $(A \setminus X) \cup \{1\}$ is a basis of
order $h' \leq h \mu$. Hence\\
$\bullet$ either $1 \in A \setminus X$, in which case we have $(A
\setminus X) = (A \setminus X) \cup \{1\}$ and then $G(A \setminus
X) = h' \leq h \mu \leq \frac{h \mu (h \mu + 3)}{2}$,\\
$\bullet$ or $1 \not\in A \setminus X$, in which case we have $(A
\setminus X) = ((A \setminus X) \cup \{1\}) \setminus \{1\}$,
implying (according to Theorem \ref{t1} for $k = 1$) that $\G(A
\setminus X) \leq \frac{h' (h' + 3)}{2} \leq \frac{h \mu (h
\mu + 3)}{2}$.\\
So, in this first case, we always have $\G(A \setminus X) \leq
\frac{h \mu (h \mu + 3)}{2}$ as required.~\vspace{1mm}\\
{\bf 2\textsuperscript{nd} case.} (if $x_n \leq 0$)\\
In this case, the elements of $X$ are all non-positive. Let $N$ be
a natural number large enough that can be written as a sum of $h$
elements of $A$; that is
\begin{equation}\label{eq32}
N = a_1 + \dots + a_t + \alpha_1 x_1 + \dots + \alpha_n x_n ,
\end{equation}
with $t , \alpha_1 , \dots , \alpha_n \in \mathbb{N}$, $a_1 ,
\dots , a_t \in A \setminus X$ and $t + \alpha_1 + \dots +
\alpha_n = h$.\\
Next, since the non-positive integer $(\alpha_1 x_1 + \dots +
\alpha_n x_n)$ is bounded below by $- (\alpha_1 + \dots +
\alpha_n) \mu = (t - h) \mu \geq t - h \mu$, then it is a sum of
$(h \mu - t)$ elements of the set $\{0 , -1\}$. It follows from
(\ref{eq32}) that $N$ is a sum of $h \mu$ elements of the set $(A
\setminus X) \cup \{0 , -1\} = (A \setminus X) \cup \{-1\}$. This
shows well (since $N$ is an arbitrary sufficiently large integer)
that the set $(A \setminus X) \cup \{-1\}$ is a basis of order
$\leq h \mu$. We finally conclude (like in the first case) that
$\G(A \setminus X) \leq \frac{h \mu (h \mu + 3)}{2}$ as
required.~\vspace{1mm}\\
{\bf 3\textsuperscript{rd} case.} (if $x_1 < 0$ and $x_n > 0$)\\
In this case, we have (from (\ref{eq30})) that $\mu = x_n - x_1$.
Let $N$ be a natural number large enough so that the number $(N +
h x_1)$ can be written as a sum of $h$ elements of $A$; that is
\begin{equation}\label{eq33}
N + h x_1 = a_1 + \dots + a_t + \alpha_1 x_1 + \dots + \alpha_n
x_n ,
\end{equation}
with $t , \alpha_1 , \dots , \alpha_n \in \mathbb{N}$, $a_1 ,
\dots , a_t \in A \setminus X$ and $t + \alpha_1 + \dots +
\alpha_n = h$.\\
From the identity
$$\alpha_1 x_1 + \dots + \alpha_n x_n - h x_1 = \alpha_2 (x_2 - x_1) + \alpha_3 (x_3 - x_1) + \dots + \alpha_n
(x_n - x_1) - t x_1 ,$$ we deduce (since $0 < x_2 - x_1 < x_3 -
x_1 < \dots < x_n - x_1 = \mu$ and $0 < - x_1 \leq x_n - x_1 - 1 =
\mu - 1$) that
$$0 < \alpha_1 x_1 + \dots + \alpha_n x_n - h x_1 \leq (\alpha_2 + \dots + \alpha_n) \mu + t (\mu - 1)
\leq h \mu - t ,$$ which implies that the integer $(\alpha_1 x_1 +
\dots + \alpha_n x_n - h x_1)$ can be written as a sum of $(h \mu
- t)$ elements of the set $\{0 , 1\}$. It follows from
(\ref{eq33}) that $N$ is a sum of $h \mu$ elements of the set $(A
\setminus X) \cup \{0 , 1\} = (A \setminus X) \cup \{1\}$. This
shows that the set $(A \setminus X) \cup \{1\}$ is a basis of
order $\leq h \mu$ and leads (as in the first case) to the desired
estimate $\G(A \setminus X) \leq \frac{h \mu (h \mu + 3)}{2}$. The
proof is
complete.\penalty-20\null\hfill$\blacksquare$\par\medbreak
\begin{rmq}\label{r1}
By using Theorem \ref{t1} of Nash for $k = 1 , 2$, we can also
establish by an elementary way (like in the above proof of Theorem
\ref{t5}) an upper bound for $\G(A \setminus X)$ in function of
$h$ and $d$. Actually, we obtain
$$\G(A \setminus X) \leq \frac{h d (h d + 1) (h d + 5)}{6} .$$
But this estimate is weaker than that of Theorem \ref{t3} and in
addition it is not linear in $d$.
\end{rmq}
{\bf Some open questions:}
\begin{enumerate}
\item[(1)] Does there exist an upper bound for $\G(A \setminus
 X)$, depending only on $h$ and $d$, which is polynomial in $h$
\underline{with degree $2$} and linear in $d$? (This asks about
the improvement of Theorem \ref{t3}).
 \item[(2)] Does there exist an upper bound for $\G(A \setminus
 X)$, depending only on $h$ and $\mu$, which is polynomial in $h$ with
 degree $2$ and \underline{linear in $\mu$}? (This asks about the improvement of Theorem
 \ref{t5}).
\end{enumerate}

\end{document}